\documentclass[11pt]{article}
\usepackage{latexsym,amssymb,amsmath}
\begin{document}
\baselineskip=20pt

\newcommand{\es}{\epsilon}
\newcommand{\ol}{\overline}
\newcommand{\la}{\langle}
\newcommand{\ra}{\rangle}
\newcommand{\psp}{\vspace{0.4cm}}
\newcommand{\pse}{\vspace{0.2cm}}

\title{\bf  Generalized  Projective  Representations
 for sl(n+1)\footnote{Research
supported by NSFC Grant 10701002}}

\author{ $\mbox{Yufeng Zhao}^1$  and  $\mbox{Xiaoping Xu}^2$\\\\\\ 1 LMAM, School of
Mathematical Sciences \\Peking University,
  Beijing, 100871, P. R. China\\\\ 2 Hua Loo-Keng Key Mathematical Laboratory\\Institute of
Mathematics, Academy of Mathematics and Systems Sciences\\ Chinese
Academy of  Sciences, Beijing, 100190, P. R. China.}
\date{ }
\maketitle
\begin {abstract}
\quad

It is well known that $n$-dimensional projective group gives rise to
a non-homogenous representation of the Lie algebra $sl(n+1)$ on the
polynomial functions of the projective space. Using Shen's mixed
product for Witt algebras (also known as Larsson functor), we
generalize the above representation of $sl(n+1)$ to a non-homogenous
representation on the tensor space of any finite-dimensional
irreducible $gl(n)$-module with the polynomial space. Moreover, the
structure of such a representation is completely determined by
employing projection operator techniques and  well-known Kostant's
characteristic identities for certain matrices with entries in
 the universal enveloping algebra. In particular, we
obtain a new one parameter family of infinite-dimensional
irreducible
 $sl(n+1)$-modules, which are in general not highest-weight type, for any given finite-dimensional irreducible
 $sl(n)$-module.
 The results could also be used to study the quantum field theory with the projective group as
the symmetry.

\end{abstract}

\section {Introduction}

A projective transformation  on $\mathbb{F}^n$ for a field
$\mathbb{F}$ is given by
$$u\mapsto \frac{Au+\vec b}{\vec c\:^t u+d}\qquad\mbox{for}\;\;u\in
\mathbb{F}^n,\eqno(1.1)$$ where all the vectors in $\mathbb{F}^n$
are in column form and
$$\left(\begin{array}{cc}A&\vec b\\ \vec c\:^t&
d\end{array}\right)\in GL(n).\eqno(1.2)$$ It is well-known that a
transformation of mapping straight lines to lines must be a
projective transformation. The group of projective transformations
is the fundamental symmetry of $n$-dimensional projective geometry.
Physically, the group  with $n=4$ and $\mathbb{F}=\mathbb{R}$
consists of all the transformations of keeping free particles
including light signals moving with constant velocities along
straight lines (e.g., cf. [GWZ1-2]). Based on the embeddings of the
poincar\'{e} group and De Sitter group into the projective group
with $n=4$ and $\mathbb{F}=\mathbb{R}$, Guo, Wu and Zhou [GWZ1-2]
proposed three kinds of special relativity.

In this paper, we give a representation-theoretic exploration on the
impact of projective transformations. Note that the Lie algebra of
$n$-dimensional projective group is spanned by the following
differential operations
$$\{\partial_{x_j},x_i\partial_{x_j},x_i\sum_{r=1}^nx_r\partial_{x_r}\mid
i,j=1,2,...,n\},\eqno(1.3)$$ which is isomorphic to the special
linear Lie algebra $sl(n+1)$. Through the above operators, we obtain
a representation of $sl(n+1)$ on the polynomial algebra
$\mathbb{F}[x_1,...,x_n]$. The non-homogeneity of (1.3) motivates us
to generalize the above representation of $sl(n+1)$ to a
non-homogenous representation on the tensor space of any
finite-dimensional irreducible $gl(n)$-module with
$\mathbb{F}[x_1,...,x_n]$ via Shen's mixed product for Witt algebras
(cf. [Sg1-3]) (also known as Larsson functor (cf. [La])). It turns
out that the structure of such generalized projective
representations can be completely determined by employing projection
operator techniques (cf. [Gm1]) and  well-known Kostant's
characteristic identities for certain matrices with entries in  the
universal enveloping algebra (cf. [K]). In particular, we obtain a
new one parameter family of infinite-dimensional irreducible
 $sl(n+1)$-modules for any given finite-dimensional irreducible
 $sl(n)$-module.

  Denote by  $\mathbb Z $ the ring of integers and by $
\mathbb N $ the additive semigroup of nonnegative integers. For any
two integers $m$ and $n$, we denote

$$\ol{m,n}=\left\{\begin{array}{ll}\{m,m+1,\cdots,n\}&\mbox{if}\;m\leq n,\\
\emptyset&\mbox{otherwise}.\end{array}\right.\eqno(1.4)$$
 Let
${\cal A}$ be a commutative associative algebra over $\mathbb{F}$.
If $\{D_{i} \ | \ i \in \overline{1,n} \}$ is a set of commuting
derivations of ${\cal A}$, the set of derivations ${\cal
W}(n)=\{\sum\limits_{i=1}^{n}a_{i}D_{i} | \ a_{i} \in {\cal A}  \}$
forms a Lie algebra via the following Lie brackets:
$$[\sum\limits_{i=1}^{n}a_{i}D_{i},\sum\limits_{i=1}^{n}b_{i}D_{i}]=
\sum\limits_{i,j=1}^{n}(a_{j}D_{j}
(b_{i})-b_{j}D_{j}(a_{i}))D_{i}.\eqno(1.5)$$

Let $E_{r,s}$ be the square matrix with 1 as its $(r,s)$-entry and 0
as the others. The general linear Lie algebra $gl(n)$ is the Lie algebra of $n\times n$
matrices over ${\mathbb F}$ with a vector space basis $\{E_{i,j}\mid
i, j \in \overline{1,n}\}$.  For any $gl(n)$-module $V$, we define
an action $\pi$ of the Lie algebra ${\cal W}(n)$ on ${\cal
A}\otimes_{\mathbb{F}}V$ by
$$\pi(\sum\limits_{i=1}^{n}a_{i}D_{i})=\sum_{i,j=1}^nD_i(a_j)\otimes E_{i,j}
+\sum\limits_{i=1}^{n} a_{i}D_{i}\otimes \mbox{Id}_V.\eqno(1.6)$$
Then $\pi$ gives a representation of ${\cal W}(n)$ (cf. [Sg1-3]) and
the functor from $gl(n)$-modules to ${\cal W}(n)$-modules was also
later known as Larsson functor (cf. [L]). The structure of the
module ${\cal A}\otimes_{\mathbb{F}}V$ was determined by Rao [R]
when ${\cal A}=\mathbb{F}[x_1^{\pm 1},...,x_n^{\pm 1}]$,
$D_i=\partial_{x_i}$ and $V$ is a finite-dimensional irreducible
$gl(n)$-module. Lin and Tan [LT] did the similar thing when ${\cal
A}$ is the algebra of quantum torus. The first author of this paper
[Z] determined the ${\cal W}(n)$-module structure of ${\cal
A}\otimes_{\mathbb{F}}V$ in the case that ${\cal A}$ is a certain
semi-group algebra and $D_i$ are locally-finite derivations in [X].
Throughout this paper, we always assume $\mbox{char}\:\mathbb{F}=0$.

 Take ${\cal
H}=\sum\limits_{i=1}^n\mathbb{F}E_{i,i}$ as a Cartan subalgebra of
$gl(n)$. Assume ${\cal A}=\mathbb{F}[x_1,...,x_n]$ and let
$D_i=\partial_{x_i}$. Embed $sl(n+1)$ into ${\cal W}(n)$ via (1.3).
Then the space ${\cal A}\otimes_{\mathbb{F}}V$ forms an
$sl(n+1)$-module with the representation $\pi|_{sl(n+1)}$, which we
call a {\it generalized projective representation}.

Characteristic identities have a long history. The first person to
exploit them was Dirac [D], who wrote down what amounts to the
characteristic identity for the Lie algebra $so(1,3)$. This
particular example is intimately connected with the problem of
describing the structure of  relativistically invariant wave
equations.
 Such identities have been shown to be  powerful tools for the analysis
of finite dimensional representations of Lie groups (cf.  [BB],
[BG], [F]). It has been shown by Kostant [K] (also cf. [Gm4]) that
the characteristic identities for semi-simple Lie algebras also hold
for infinite dimensional representations. Moreover, one may
construct projection operators analogous to the projection operators
of Green [G], and Bracken and Green [BG] in finite dimensions. We
refer
  [Gh], [Gm1-Gm3], [Ha],  [L], [LG], [M], [OCC] and [O] for the
 other works on the identities. Using
projection operator techniques and Kostant's characteristic
identities, we prove:

\vspace{0.2cm}
 {\bf
Main Theorem}. {\it Let $V$ be a finite-dimensional irreducible
$gl(n)$-module with highest weight  $\mu$. We have the following
conclusions:

(i) The space ${\cal A}\otimes_{\mathbb{F}}V$ is an irreducible
$sl(n+1)$-module if and only if
$$\mu(E_{1,1})+\sum\limits_{j=1}^{n}\mu(E_{j,j})\not\in
-\mathbb{N}\bigcup\ol{2,1+\mu(E_{1,1}-E_{2,2})}\eqno(1.7)$$ and
$$\mu(E_{i,i})+\sum\limits_{j=1}^{n}\mu(E_{j,j})-i\not\in\ol{1,
\mu(E_{i,i}-E_{i+1,i+1})}\qquad\mbox{for}\;i\in\ol{2,n-1}.\eqno(1.8)$$

(ii) If one of the conditions in (1.7) and (1.8) fails, then both
$U(sl(n+1))(1\otimes V)$ and $({\cal
A}\otimes_{\mathbb{F}}V)/(U(sl(n+1))(1\otimes V))$ are irreducible $sl(n+1)$-modules .}\vspace{0.2cm}

Note that all $\mu(E_{i,i}-E_{i+1,i+1})$ are nonnegative integers,
which determine the corresponding $sl(n)$-module $V$ uniquely.
Moreover, the identity matrix in $gl(n)$ are allowed to be any
constant map. The above theorem says that given a finite-dimensional
irreducible $sl(n)$-module, we can construct a new one-parameter
family of explicit irreducible $sl(n+1)$-modules via its projective
representation and Shen's mixed product.

 A quantum field is an operator value function on a certain Hilbert space,
which is often a direct sum of  infinite-dimensional irreducible
modules of a certain Lie algebra (group). The Lie algebra of
two-dimensional conformal group is exactly the Virasoro algebra,
which is infinite-dimensional. The minimal models of two-dimensional
conformal field theory were constructed from direct sums of certain
infinite-dimensional irreducible modules of the Virasoro algebra,
where a distinguished module gives rise to a vertex operator
algebra. When $n>2$, the $n$-dimension conformal group is
finite-dimensional, whose Lie algebra is exactly isomorphic to
$so(n,2)$. It is still unknown what should a higher-dimensional
conformal field theory be. Part of reason is that we lack of enough
knowledge on the infinite-dimensional irreducible $so(n,2)$-modules
that are compatible to the natural conformal representation of
$so(n,2)$. This motivates us to study explicit infinite-dimensional
irreducible modules of finite-dimensional simple Lie algebras by
using non-homogeneous polynomial representations and Shen's mixed
product for Witt algebras. This paper is the first work in this
direction.

As we mentioned earlier, projective groups are important groups in
physics. In comparison with the minimal models of two-dimensional
conformal field theory, the underlying module of the projective
representation of $sl(n+1)$ should be the distinguished  module in
the possible quantum field theory with the projective group as the
symmetry. The other modules would make the theory more substantial.

 The paper is organized as follows. In Section 2, we
slightly generalize Kostant's characteristic identities and recall
some facts about projection operators based on Kostant's work [K]
and Gould's works [Gm1, Gm4]. In Section 3, we prove (i) and (ii) in
the theorem. \vspace{0.2cm}

{\bf Acknowledgement:} We would like to thank Professor Han-Ying Guo
for his interesting talk that motivates this work.

\section {Characteristic  Identities and Projection Operators}

\quad \quad In order to keep the paper self-contained, we will first prove
certain characteristic identities
  for $gl(n)$, which will be used to study  the irreducibility of the generalized projective representations for $sl(n+1)$ . Then we will recall some facts about the projection operators for $gl(n)$, although part of
 the results in this section has appeared in  [Gm1, Gm4], [OCC] and [K].

\subsection{Some Standard Facts for $gl(n)$ and $sl(n)$ }

\quad \quad Recall that $gl(n)$ is the Lie algebra of $n\times n$
matrices over ${\mathbb F}$ with a  basis $\{E_{i,j}\mid i, j \in
\overline{1,n}\}$ and the Lie bracket:
$$[E_{i,j},E_{k,l}]=\delta_{k,j}E_{i,l}-\delta_{i,l}E_{k,j}\qquad\mbox{for}\;\;
i,j,k,l \in \overline{1,n}.\eqno(2.1)$$  Note that $gl(n)$ is
reductive  with the following decomposition of ideals
$gl(n)=sl(n)\oplus {\mathbb F}\mbox{I}$, where $sl(n)$ is the
special linear Lie algebra of  matrices with zero trace and
$\mbox{I}=\sum\limits_{i=1}^n E_{i,i}$ is the identity matrix, which
is central.

 Take ${\cal
H}=\sum\limits_{i=1}^n\mathbb{F}E_{i,i}$ as  a Cartan subalgebra of
$gl(n)$. If $\lambda \in {\cal{H}}^{*}$ is a weight of $gl(n)$, we
identify $\lambda$ with the $n$-tuple
$\lambda=(\lambda_{1},\cdots,\lambda_{n})$, where
$\lambda_{i}=\lambda(E_{i,i})$. For $\lambda, \mu  \in
{\cal{H}}^{*}$, we define $$(\lambda,
\mu)=\sum\limits_{i=1}^{n}\lambda_{i}\mu_{i}.\eqno(2.2)$$

 Denote by $\varepsilon_i$ the weight with 1 as its $i$th coordinate
 and 0 as the others, i.e. $$\varepsilon_{i}=(0,\cdots, 0, \stackrel{i}{1}, 0, \cdots, 0).\eqno(2.3)$$
 The set $\Phi^+=\{\varepsilon_{i}-\varepsilon_{j}\mid 1\leq i <
 j\leq n\}$ forms a set of positive roots of $gl(n)$. In this case, the
half-sum of the positive roots is given
by$$\delta=\frac{1}{2}\sum\limits_{i <
j}(\varepsilon_{i}-\varepsilon_{j})=\frac{1}{2}\sum\limits_{i=1}^{n}
(n+1-2i)\varepsilon_{i}.\eqno(2.4)$$ Moreover, there exists a
one-to-one correspondence  between the set of finite-dimensional
irreducible $gl(n)$-modules and the set of $n$ tuples
$\lambda=(\lambda_{1},\cdots,\lambda_{n})$ such that
$\lambda_{i}-\lambda_{i+1} \in {\mathbb{N} } \  \mbox{for} \ i \in
\overline{1,n-1}$. Such an $n$ tuple $\lambda$ is called the {\it
highest weight} of the corresponding module which we  denote as
$V(\lambda)$.\vspace{0.2cm}

 Let $U$ denote the universal enveloping
algebra of $gl(n)$ and let $Z$ be the center of $U$. Set
$$\sigma_{1}=\mbox{I}, \ \sigma_r=\sum\limits_{i_1,..,i_r=1}^n
E_{i_1,i_2}E_{i_2,i_3}\cdots
E_{i_{r-1},i_r}E_{i_r,i_1}\qquad\mbox{for}\;r\in\ol{2,n}.\eqno(2.5)$$
Then the center
$$Z=\mathbb{F}[\sigma_{1},\cdots,\sigma_{n}].\eqno(2.6)$$
The subspace
${\cal H}_0={\cal H}\bigcap sl(n)$ is a Cartan
subalgebra of $sl(n)$ with the standard basis
$\{\alpha_{i}^\vee=E_{i,i}-E_{i+1,i+1}\mid i \in
\overline{1,n-1}\}.$ Denote the dual vector space of ${\cal
H}_{0}$ by ${\cal H}_{0}^*$ and let
$\omega_{1},\omega_{2},\cdots,\omega_{n-1}$ be the fundamental
integral dominant weights in ${\cal H}_{0}^*$ defined by
$\omega_{i}(\alpha_{j}^\vee)=\delta_{i,j}.$\vspace{0.3cm}

Let $V(\psi)$ be the  finite dimensional irreducible $sl(n)$-module
with the highest weight $\psi$. We can make $V(\psi)$ as a
$gl(n)$-module $V(\psi,b)$ by letting the central element $\mbox{I}$ act
as the scalar map $b\mbox{Id}_{V(\psi)}$.

For $\vec a=(a_1,...,a_{n-1})\in\mathbb{N}^{n-1}$ and
$0<k\in\mathbb{Z}$, we denote \begin{eqnarray*}\hspace{1.4cm}I(\vec
a,k)&=&\{(a_1+c_1-c_2,a_2+c_2-c_3,...,a_{n-1}+c_{n-1}-c_n)\mid
c_i\in\mathbb{N}\\ & &\mbox{such
that}\;\sum_{i=1}^nc_i=k\;\mbox{and}\;c_{s+1}\leq
a_s\;\mbox{for}\;s\in\ol{1,n-1}\}.\hspace{3cm}(2.7)\end{eqnarray*}
Moreover, we set
$$\omega_{\vec a}=\sum_{i=1}^{n-1}a_i\omega_i\qquad\mbox{for}\;
\vec a\in\mathbb{N}^{n-1}.\eqno(2.8)$$\vspace{0.1cm}

{\bf Lemma 2.1.1} (e.g., cf. Proposition 15.25 in  [FH]) {\it \quad
For any $\vec a\in\mathbb{N}^{n-1}$, the tensor product of
$sl(n)$-module $V(\omega_{\vec a})$ with
 $V(k\omega_{1})$ decomposes into a direct sum:
$$V(\omega_{\vec a})\otimes_{\mathbb{F}}
V(k\omega_{1})=\bigoplus_{\vec b\in I(\vec a,k)} V(\omega_{\vec
b}).\eqno(2.9)$$}
\vspace{0.1cm}

{\bf Lemma 2.1.2} {\it \quad  Let  $\Pi$ be the weight set of $gl(n)$-module $V(\psi,b)$. Assume $\psi=\sum\limits_{i=1}^{n-1}a_i\omega_i$, $\nu=(\nu_{1},\cdots,\nu_{n})\in \Pi$ and $(\nu_{1}-\nu_{2},\cdots,\nu_{n-1}-\nu_{n})=\psi-\sum\limits_{i=1}^{n-1}k_i\alpha_i$. Then
$$\nu_{1}=\sum\limits_{i=1}^{n-1}a_i+\frac{b-\sum\limits_{i=1}^{n-1}ia_{i}}{n}-k_{1},
 \ \nu_{n}=\frac{b-\sum\limits_{i=1}^{n-1}ia_{i}}{n}+k_{n-1}, $$
 $$ \nu_{j}=\sum\limits_{i=j}^{n-1}a_i+\frac{b-\sum\limits_{i=1}^{n-1}ia_{i}}{n}+k_{j-1}-k_{j},\;\;j \in \overline{2,n-1}.\eqno(2.10) $$  }\vspace{0.1cm}

Denote
$$\underline{k}=(k_1,k_2,...,k_n)\in \mathbb{N}^n,\quad
|\underline{k}|=\sum_{i=1}^nk_i, $$$$I(\mu,j)=\{\underline{c}=(c_{1},\cdots,c_{n}) \  | \  c_i \in\mathbb{N} \ , \ |\underline{c}|=j , \ c_{s+1}\leq
\mu_{s}-\mu_{s+1} \;\mbox{for}\;s\in\ol{1,n-1}\}. \eqno(2.11)$$

It is easy to deduce the following lemma from  the above two
lemmas.\vspace{0.3cm}

{\bf Lemma 2.1.3} {\it \quad  The tensor product of $gl(n)$-module $V(\mu)$ with
 $V(k\varepsilon_{1})$ decomposes into a direct sum:
$$V(\mu)\otimes_{\mathbb{F}}
V(k\varepsilon_{1})=\bigoplus_{\underline{c}\in I(\mu,k)} V(\mu+\underline{c}).\eqno(2.12)$$}\vspace{0.3cm}

\subsection{Characteristic Identities and Projection Operators for $gl(n)$}

\quad Now we will introduce the characteristic identities for $gl(n)$.
We know that the universal enveloping algebra $U$ of $gl(n)$ can be imbedded
into $U\otimes U$ by the associative algebra homomorphism $d: U
\rightarrow U\otimes U$ determined  by $$d(u)=u\otimes 1 +1 \otimes
u \qquad \mbox{ for} \  u\in gl(n).\eqno(2.13)$$
 Let $V(\lambda)$ be a fixed finite-dimensional $gl(n)$-module  with highest weight $\lambda$
  and let $\pi_{\lambda}$ be the corresponding
representation. Kostant [K] considered the map
$$\partial: U \rightarrow (\mbox{End}\:V(\lambda))\otimes_{\mathbb{F}} U; $$$$\ u
\mapsto 1\otimes u+\pi_{\lambda}(u)\otimes 1\eqno(2.14)$$ for $u \in
gl(n)$ and extended $\partial$ to an associative algebra
homomorphism from $U$ to
$(\mbox{End}\:V(\lambda))\otimes_{\mathbb{F}} U$. More generally, if
$d(u)=\sum\limits_{r}u_{r}\otimes v_{r}$, we have $\partial(u)=
\sum\limits_{r}\pi_{\lambda}(u_{r})\otimes v_{r}$. For $z \in Z$, we
denote  $$
\tilde{z}=-\frac{1}{2}[\partial(z)-\pi_{\lambda}(z)\otimes
1-1\otimes z],\eqno(2.15)$$which may be viewed as an $m \times m$
($m=\mbox{dim} V(\lambda)$)  matrix with entries in $U$.

Denote by $\chi_\zeta$ the central character of a highest weight
$gl(n)$-module with highest weight $\zeta$. Suppose now that $W$ is
another $gl(n)$-module  admitting the central character $\chi_{\mu}$
and $\pi_{\mu}$ is the corresponding representation. We extend
$\pi_{\mu}$ to an algebra homomorphism
$$ \tilde{\pi_{\mu}}: (\mbox{End}\:V(\lambda)) \otimes_{\mathbb{F}} U \rightarrow
(\mbox{End}\:V(\lambda)) \otimes_{\mathbb{F}} (\mbox{End}\:W); $$$$\
\sum\limits_{i}\rho_{i}\otimes u_{i} \mapsto
\sum\limits_{i}\rho_{i}\otimes \pi_{\mu}(u_{i}),\eqno(2.16)
$$where $\rho_{i} \in \mbox{End}\:V(\lambda)$ and $u_{i} \in U$. In particular, we have
$$\tilde{\pi_{\mu}}(\tilde{z})=
-\frac{1}{2}[(\pi_{\lambda}\otimes \pi_{\mu})
(z)-\pi_{\lambda}(z)\otimes 1-1\otimes \pi_{\mu}(z)]\eqno(2.17)$$
for $z \in Z$. Clearly, $\tilde{\pi_{\mu}}(\tilde{z})$ is a linear
operator on $V(\lambda) \otimes_{\mathbb{F}} W$ which may be viewed
as an $m \times m$  matrix with entries from $\mbox{End}\:W$ under a
 basis of $V(\lambda)$. For $\nu\in {\cal H}^\ast$, we define
$$f_\nu=-\frac{1}{2}(\chi_{\mu+\nu}-\chi_\lambda-\chi_\mu).\eqno(2.18)$$
Denote by $\Pi_\lambda$ the weight set of $V(\lambda)$.
\vspace{0.3cm}

{\bf Lemma 2.2.1} (cf.  [K], [G], [OCC] ) \quad {\it On the space
$W$, the matrix $\tilde{z}$ satisfies the following characteristic
identity:
$$\prod\limits_{\nu\in\Pi_\lambda}(\tilde{z}-f_\nu(z))=0
\qquad\mbox{for}\;\;z\in Z.\eqno(2.19)$$ }

By varying the module $V(\lambda)$ and the central element $z$, we
obtain a series of characteristic  identities. In particular, the
following characteristic  identity will be used in the proof of the
main theorem:\vspace{0.3cm}

 {\bf Corollary  2.2.2}\quad {\it
 Take $V(\lambda)$ to be the dual module of
 $gl(n)$-module $V(2,1,\cdots,1)$.  Then the matrix $\tilde{\sigma_{2}}$ satisfies the
following characteristic identity on $W$:
$$\prod\limits_{i=1}^{n}(\tilde{\sigma_{2}}-m_i)=0,
\ \
m_i=\frac{1}{2}(\lambda,\lambda+2\delta)-\frac{1}{2}(\lambda_{i},
\lambda_{i}+2(\mu+\delta)),\eqno(2.20)$$
where $\lambda_{i}=(-1,\cdots, \stackrel{i}{-2},\cdots, -1), \ i \in
\overline{1,n}$.

}\vspace{0.3cm}

{\it Proof}\quad Note that the module $V(2,1,\cdots,1)$ has a basis
$\{e_{1},\cdots,e_{n}\}$ such that
$$E_{i,i}(e_{k})=(1+\delta_{i,k})e_{k}, \
E_{i,j}(e_{k})=\delta_{j,k}e_{i} \ (i \neq j).\eqno(2.21)$$
 Let $\pi^{*}$ be
the dual module of $gl(n,{\mathbb F})$-module $V(2,1,\cdots,1)$.
Then $$\pi^{*}(E_{i,i})=-(E_{i,i}+\mbox{I}), \
\pi^{*}(E_{i,j})=-E_{j,i} \ (i \neq j).\eqno(2.22)$$ Obviously, the
representation $\pi^{*}$ is $n$-dimensional and all its weights
$$\{\lambda_{1}=(-2,-1,\cdots,-1),..., \lambda_{i}=(-1,\cdots,
\stackrel{i}{-2},\cdots, -1), \cdots,
\lambda_{n}=(-1,\cdots,-1,-2)\} \eqno(2.23)$$occur with multiplicity
one. The matrix $\tilde{\sigma_{2}}$ in this case is given by
$$\tilde{\sigma_{2}}=-\sum\limits_{i,j}^{n}\pi^{*}(E_{j,i})E_{i,j}.\eqno(2.24)$$
This is the matrix
$$\begin{array}{rcl}
\tilde{\sigma_{2}}&=& \left[\begin{array}{cccc}
\mbox{I}+E_{1,1}&E_{1,2}&\cdots&E_{1,n}\\
E_{2,1}&\mbox{I}+E_{2,2}&\cdots&E_{2,n}\\
\vdots&\vdots&\cdots&\vdots\\
E_{n,1}&E_{n,2}&\cdots&\mbox{I}+E_{n,n}
\end{array}
\right].
\end{array} \eqno(2.25)$$
By Lemma 2.2.1, the operator $\tilde{\sigma_{2}}$ satisfies the
characteristic identity on $W$:
$$\prod\limits_{i=1}^{n}(\tilde{\sigma_{2}}-m_i)=0,
\ \
m_i=\frac{1}{2}(\lambda,\lambda+2\delta)-\frac{1}{2}(\lambda_{i},
\lambda_{i}+2(\mu+\delta)),
 \ i \in
\overline{1,n}.\eqno(2.26)$$\hspace{1cm} $\Box$\vspace{0.3cm}

In the rest of this section, we will recall some facts about
projection operators appeared in [Gm2].

 Take $gl(n)$-module $V(\lambda)=V(\varepsilon_{1})^{*}$ (resp. $V(\varepsilon_{1})$). Then the corresponding
 matrices
 $\tilde{\sigma_{2}}$ are
 $$M=
(E_{i,j})_{i,j=1}^{n}, \ \tilde{M}=-M^{T} \in U(gl(n)),\eqno(2.27)$$
respectively. On the space $W$, the matrices $M$ and $\tilde{M}$
satisfy the following characteristic identities:
$$\prod\limits_{i=1}^{n}(M-d_i)=0, \ d_{i}=\mu_{i}+n-i;\eqno(2.28)$$
$$\prod\limits_{i=1}^{n}(\tilde{M}-\tilde{d}_i)=0, \ \tilde{d}_{i}=n-1-d_{i}, \ i \in \overline{1,n}.\eqno(2.29) $$
For $r \in \overline{1,n}$,
$$P_{r}=\prod\limits_{r\neq l \in \overline{1,n}}(\frac{M-d_{l}}{d_{r}-d_{l}}),\qquad \tilde{P}_{r}=\prod\limits_{r\neq l \in \overline{1,n}}(\frac{\tilde{M}-\tilde{d}_{l}}{\tilde{d}_{r}-\tilde{d}_{l}})\eqno(2.30)$$are called {\it  projection  operators}, which  project the tensor product space
$V(\varepsilon_{1})^{*}{\otimes}_{\mathbb{F}} W$ (resp. $V(\varepsilon_{1})\otimes_{\mathbb{F}} W$ ) onto the irreducible
 module $W_{r}=P_{r}(V(\varepsilon_{1})^{*}\otimes_{\mathbb{F}} W)$ (resp. $\tilde{W}_{r}=\tilde{P}_{r}(V(\varepsilon_{1})\otimes_{\mathbb{F}} W)$)
 with central character $\chi_{\nu-\varepsilon_{r}}$ (resp. $\chi_{\nu+\varepsilon_{r}}$).

\section{Generalized Projective Representations}

\quad \quad In this section, we will   give the
  detailed
  construction of generalized projective representations
 for  the special linear Lie
 algebra $sl(n+1)$ and study their irreducibility.

  \subsection{Construction of the Representations}

\qquad Let ${\cal A}=\mathbb{F}[x_{1},\cdots,x_{n}]$. The operator
$p_{i}=x_{i}\sum\limits_{i=1}^{n}x_{i}\partial_{x_{i}}$ is called
{\it pseudo-translation operator} on ${\cal A}$ in physics. Note
that the Lie algebra of $n$-dimensional projective group
$$L_{n+1}=\sum\limits_{i,j=1}^{n}\mathbb{F}x_{i}\partial_{x_{j}}+
\sum\limits_{i=1}^{n}\mathbb{F}\partial_{
x_{i}}+\sum\limits_{i=1}^{n}\mathbb{F}p_{i},\eqno(3.1)$$ forms a Lie
subalgebra of Witt algebra
$${\cal W}(n)=\{\sum\limits_{i=1}^{n}f_{i}\partial_{ x_{i}} \ | \
f_{i} \in {\cal A} \} .\eqno(3.2) $$
 Moreover, we have
the following Lie brackets:
$$[\partial_{x_{j}},p_{i}]=\left\{\begin{array}{llll}
x_{i}\partial_{x_{j}},  & i \neq j,\\
\sum\limits_{i=1}^{n}x_{i}\partial_{x_{i}}+x_{i}\partial_{x_{i}}, &
i
=j,\end{array}\right.\eqno(3.3)$$$$[x_{i}\partial_{x_{j}},p_{k}]=\delta_{j,k}{p}_{i},
\
[x_{i}\partial_{x_{j}},\partial_{x_{k}}]=-\delta_{i,k}\partial_{x_{j}},
\
[x_{i}\partial_{x_{j}},x_{k}\partial_{x_{l}}]=\delta_{j,k}x_{i}\partial_{x_{l}}-
\delta_{i,l}x_{k}\partial_{x_{j}}.\eqno(3.4)$$To abbreviate,  we
denote
$$P=\sum\limits_{i=1}^{n}\mathbb{F}p_{i}, \
S=\sum\limits_{i=1}^{n}\mathbb{F}\partial_{ x_{i}}, \
\overline{L}_{n}=\sum\limits_{i,j=1}^{n}\mathbb{F}x_{i}\partial_{x_{j}}, \ \overline{L}_{n}'=[\overline{L}_{n},\overline{L}_{n}].\eqno(3.5)$$
Then
$$L_{n+1}=P \oplus S \oplus \overline{L}_{n}, \  [P,P]=\{0\}, \
[S,S]=\{0\}.
 \eqno(3.6)$$
Moreover, $\overline{L}_{n}$ (resp. $\overline{L}_{n}'$) is
isomorphic to $gl(n)$ (resp. $sl(n)$).  It is easy to
verify:\vspace{0.3cm}

 {\bf Lemma 3.1.1 } {\it The special linear Lie algebra $sl(n+1)$ is isomorphic
to $L_{n+1}$ with the following identification of Chevalley
generators:
$$h_{i}=x_{i}\partial_{x_{i}}-x_{i+1}\partial_{x_{i+1}},  \qquad
h_{n}=\sum\limits_{i=1}^{n}x_{i}\partial_{x_{i}}+x_{n}\partial_{x_{n}},\eqno(3.7)$$
$$e_{i}=x_{i}\partial_{x_{i+1}}, \;\; f_{i}=x_{i+1}\partial_{x_{i}},
\qquad e_{n}=p_{n}, \;\; f_{n}=-\partial_{x_{n}},\eqno(3.8)$$ for $i
\in \overline{1,n-1}$. }\vspace{0.1cm}

For any finite dimensional $gl(n)$-module $V$ with highest weight
$\mu=(\mu_{1},\cdots,\mu_{n}) $, we define an action $\pi$ of Witt
Lie algebra ${\cal W}(n)$ on ${\cal A}\otimes_{\mathbb{F}}V$ by
$$\pi(\sum\limits_{i=1}^{n}a_{i}D_{i})=\sum_{i,j=1}^nD_i(a_j)\otimes E_{i,j}
+\sum\limits_{i=1}^{n} a_{i}D_{i}\otimes \mbox{Id}_V.\eqno(3.9)$$
Embed $sl(n+1)$ into ${\cal W}(n)$ by (3.7) and (3.8). The space
${\cal A}\otimes_{\mathbb{F}}V$ forms an $sl(n+1)$-module with the
representation $\pi|_{sl(n+1)}$, which we call a {\it generalized
projective representation} of $sl(n+1)$.\vspace{0.3cm}

It follows from (3.7)-(3.9) that the explicit $L_{n+1} $-module ${\cal A}\otimes_{\mathbb{F}}V$
structure is given by:

$$(x_{i}\partial_{x_{j}}).(f\otimes v)=(x_{i}\partial_{x_{j}})f\otimes
v +f\otimes E_{i,j}.v,
 \eqno(3.10)$$$$\partial_{x_{i}}.(f\otimes v)=\partial_{x_{i}}(f)\otimes v,\eqno(3.11)$$$$
p_{i}.(f\otimes v)=p_{i}(f)\otimes v+x_{i}f\otimes
\sum\limits_{i=1}^n E_{i,i}.v+\sum\limits_{j=1}^n x_{j}f\otimes
E_{i,j}.v \eqno(3.12)$$ for $i,j\in\ol{1,n}$, where $f \in {\cal A}$
and  $v \in V$.

\subsection{Irreducibility Criteria for Generalized Projective Representations}

 \qquad In this section,  we will give two irreducibility criteria for  $L_{n+1}$-module  ${\cal A}\otimes_{\mathbb{F}}V$ (cf. Proposition 3.2.3 and Proposition 3.2.5).\vspace{0.3cm}

{\bf Lemma 3.2.1} \quad {\it The vector space $U(P)(1\otimes_{\mathbb{F}} V)$ is an
 irreducible $L_{n+1}$-submodule of ${\cal A}\otimes_{\mathbb{F}}V$,
where $U(P)$ denote the universal enveloping algebra of abelian Lie
algebra $P$.}\vspace{0.3cm}

 {\it Proof} \quad Denote
$$x^{\underline{c}}=x_{1}^{c_1}x_{2}^{c_2}\cdots x_{n}^{c_n}  \qquad \mbox{for} \ \underline{c}=(c_{1},\cdots,c_{n}) \in \mathbb{N}^{n}.\eqno(3.13)$$
Recall that $\Pi_\mu$ denotes the weight set of $V$. By (3.10), we
have
$$(x_{i}\partial_{x_{i}}).(x^{\underline{c}}\otimes v_{\nu})=(c_{i}+\nu_{i})
x^{\underline{c}}\otimes v_{\nu}, \eqno(3.14)$$
$$(\sum\limits_{i=1}^{n}x_{i}\partial_{x_{i}}).(x^{\underline{c}}\otimes
v_{\nu})=(|\underline{c}|+\sum\limits_{i=1}^{n}\nu_{i}) x^{\underline{c}}\otimes
v_{\nu},
 \eqno(3.15)$$ for any $
x^{\underline{c}} \in {\cal A}$ and $\nu \in \Pi_\mu$.

Assume that
$\{v_{1},\cdots,v_{\ell}\}$ is a basis of $V$. For $1 \leq k \in
\mathbb{N}$, we set
$$ U(P)(1\otimes_{\mathbb{F}}
V)_{{\la}k\ra}=\mbox{Span}_{\mathbb{F}}\{ p_{i_{1}}p_{i_{2}}\cdots
p_{i_{k}}(1\otimes v_{j}) \ | \ 1 \leq i_{1}\leq i_{2}\leq \cdots
\leq i_{k}\leq n, j \in \overline{1,d} \}, \eqno(3.16)$$
$$({\cal A}\otimes_{\mathbb{F}}V)_{{\la}k\ra}=\mbox{Span}_{\mathbb{F}}\{x^{\underline{c}}\otimes v_{j} \ | \
x^{\underline{c}} \in {\cal A}, \ |\underline{c}|=k, \ j \in \overline{1,l}\}.
\eqno(3.17)$$
 Then (3.14) implies
that
 $$U(P)(1\otimes_{\mathbb{F}} V)=\bigoplus\limits_{k \in \mathbb{N}}(U(P)(1\otimes
V))_{{\la}k\ra}.\eqno(3.18)$$Denote
$$\triangle_{i,j}^{k}=\left\{\begin{array}{llll}
x_{j}\partial_{x_{i}}  & \mbox{if}\;i \neq j,\\
\sum\limits_{i=1}^{n}x_{i}\partial_{x_{i}}+x_{i}\partial_{x_{i}}+k-1
& \mbox{if}\;i =j,\end{array}\right.\eqno(3.19)$$By induction, we can easily verify the following fomula:
$$\partial_{x_{i}}.p_{i_{1}}p_{i_{2}}\cdots p_{i_{k}}(1\otimes
v_{j})=\sum\limits_{s=1}^{n}p_{i_{1}}p_{i_{2}}\cdots\hat{p}_{i_{s}}\cdots
p_{i_{k}}\triangle_{i,i_{s}}^{k}(1\otimes v_{j}),\eqno(3.20)$$
\begin{eqnarray*}\hspace{1cm}& &
(x_{i}\partial_{x_{l}}).p_{i_{1}}p_{i_{2}}\cdots p_{i_{k}}(1\otimes
v_{j})\\
&=&\sum\limits_{s=1}^{n}\delta_{l,i_{s}}p_{i_{1}}p_{i_{2}}\cdots\hat{p}_{i_{s}}\cdots
p_{i_{k}}p_{i}(1\otimes
v_{j})+p_{i_{1}}p_{i_{2}}\cdots{\cal}{p}_{i_{k}}(1\otimes
E_{i,l}.v_{j}),\hspace{2.6cm}(3.21)\end{eqnarray*}
$${p}_{i}.{\cal}{p}_{i_{1}}{p}_{i_{2}}\cdots p_{i_{k}}(1\otimes
v_{j})={p}_{i}p_{i_{1}}p_{i_{2}}\cdots {p}_{i_{k}}(1\otimes
v_{j}).\eqno(3.22)$$  So (3.20)-(3.22) imply that
$U(P)(1\otimes_{\mathbb{F}} V)$ is an $L_{n+1}$-submodule of  ${\cal
A}\otimes_{\mathbb{F}}V$. Furthermore, it is easy to verify that
$$\bigcap\limits_{i=1}^{n} \mbox{Ker} \ \partial_{x_{i}}|_{({\cal
A}\otimes_{\mathbb{F}}V)_{{\la}k\ra}}=\{0\} \  \mbox{for \ any } \ 1
\leq k \in \mathbb{N}. \eqno(3.23)
$$Thus any non-trivial submodule of $U(P)(1\otimes_{\mathbb{F}} V)$ must contain $1\otimes V$.
The irreducibility of  $U(P)(1\otimes_{\mathbb{F}} V)$ follows.
\hspace{1cm}
$\Box$\vspace{0.3cm}

In the rest of this section,  we will investigate the condition for
${\cal A}\otimes_{\mathbb{F}}V=U(P)(1\otimes_{\mathbb{F}}
V)$.\vspace{0.3cm}

 It is obvious that
${\cal A}\otimes_{\mathbb{F}}V=U(P)(1\otimes_{\mathbb{F}} V)$ if and
only if $({\cal
A}\otimes_{\mathbb{F}}V)_{{\la}k\ra}=(U(P)(1\otimes_{\mathbb{F}}
V))_{{\la}k\ra}$ holds for any $   k \in \mathbb{N}$.\vspace{0.3cm}

 Suppose that $B_{i}$ is an ordered basis for $({\cal
   A}\otimes_{\mathbb{F}}V)_{{\la}i\ra}$. Denote by
${P}_{i+1,i}^{j}$ the matrix of the linear map ${p}_{j}|_{({\cal
A}\otimes_{\mathbb{F}}V)_{{\la}i\ra}}$ (cf. (3.12))
  with respect to the bases $B_{i}$ and $B_{i+1}$. For any
  $0<j\in\mathbb{Z}$, we denote
  $$\Gamma_j=\{\hat i=(i_1,i_2,...,i_j)\mid i_s\in\ol{1,n};i_1\leq
  i_2\leq\cdots\leq i_j\}.\eqno(3.24)$$
Set
$$P_{\hat
i}={P}_{j,j-1}^{i_{1}}{P}_{j-1,j-2}^{i_{2}}\cdots{P}_{1,0}^{i_{j}}.\eqno(3.25)$$
We order
$$\Gamma_j=\{\hat k^1,\hat k^2,...,\hat k^{\ell_j}\}\eqno(3.26)$$
lexically. In particular,
$$\hat k^1=(1,1,...,1),\qquad \hat
k^{\ell_j}=(n,n,...,n).\eqno(3.27)$$
  Then $(U(P)(1\otimes_{\mathbb{F}} V))_{{\la}j\ra}$ is, as a vector
  space,
isomorphic to the column space of the $\ell\ell_j\times \ell\ell_j$
matrix
$$
M_{j}= [P_{\hat k^1},P_{\hat k^2},...,P_{\hat k^{\ell_j}}]
\eqno(3.28)$$ (recall that $\ell=\dim V$).
 Denote
$$I_{1}=\{1\}\bigcup\{i \in \overline{2,n} \ | \ \mu_{i-1}-\mu_{i}\geq 1\}.
\eqno(3.29)$$ By means of the characteristic identity  in Corollary
2.2.2, we get the following necessary but not sufficient condition
for the irreducibility of  $L_{n+1}$-module  ${\cal
A}\otimes_{\mathbb{F}}V$:\vspace{0.3cm}

 {\bf Lemma 3.2.2 } {\it   \quad For any $1 \leq k \in \mathbb{N}$. If $\mu_{i}+|\mu|-i+s \neq 0$ for any $ s \in
\overline{1,k}, \ i \in I_{1}$, then  $L_{n+1}$-module  ${\cal A}\otimes_{\mathbb{F}}V$ is irreducible.
 }\vspace{0.3cm}

 {\it Proof} \quad   The Lemma is based on the following result:\vspace{0.3cm}

 {\it Claim. For any $1 \leq k \in \mathbb{N}$. If $\mu_{i}+|\mu|-i+s \neq 0,  \ \forall \ s \in
\overline{1,k}, \ i \in I_{1}$, then
$({\cal A}\otimes_{\mathbb{F}}V)_{{\la}k\ra}=(U(P)(1\otimes_{\mathbb{F}}
V))_{{\la}k\ra}$ holds. }\vspace{0.3cm}

 We will prove this claim by induction on $k$.
For $k=1$, $(U(P)(1\otimes_{\mathbb{F}} V))_{{\la}1\ra}$ is
isomorphic as vector space to the column space of $\ell n\times \ell
n$ matrix
$$\begin{array}{rcl}
&& [{P}_{1,0}^{1},{P}_{1,0}^{2},\cdots,{P}_{1,0}^{n} ]
\end{array} \mbox{with} \ \begin{array}{rcl}
{P}_{1,0}^{k}&=& \left[\begin{array}{c} E_{k,1}|_V\\E_{k,2}|_V\\\vdots\\
(\mbox{I}+E_{k,k})|_V\\\vdots\\E_{k,n}|_V
\end{array}
\right],
\end{array}\eqno(3.30)$$
which is exactly $\tilde{\sigma_{2}}|_V$ (cf. (2.25)).
By Corollary
2.2.2, the matrix $\tilde{\sigma_{2}}|_V$ is diagonalizable and it has
full rank if and only if all its eigenvalues are not zero, i.e.
$$m_{i}=\frac{1}{2}(\lambda,\lambda+2\delta)-\frac{1}{2}
 (\lambda_{i},\lambda_{i}+2(\mu+\delta))
 =\mu_{i}+ |\mu|-i+1\neq
 0\eqno(3.31)$$ for any $i \in I_{1}$.

Now suppose that the lemma holds for $k=\iota-1$. Assume $k=\iota$.
Note that the eigenvalues of the matrix
$$ \left[\begin{array}{cccc}
(\mbox{I}+E_{1,1}+s-1)|_V&E_{1,2}|_V&\cdots&E_{1,n}|_V\\
E_{2,1}|_V&(\mbox{I}+E_{2,2}+s-1)|_V&\cdots&E_{2,n}|_V\\
\vdots&\vdots&\cdots&\vdots\\
E_{n,1}|_V&E_{n,2}|_V&\cdots&(\mbox{I}+E_{n,n}+s-1)|_V
\end{array}
\right]  \eqno(3.32)$$  are  $\mu_{i}+\sum\limits_{ j
=1}^{n}\mu_{j}-i+s  \ ( i \in I_{1})$ by (3.31). Thus it is
invertible if and only if $\mu_{i}+|\mu|-i+s \neq 0,  \ \forall \ i
\in I_{1}$.\vspace{0.3cm}

 Observe that
$({\cal A}\otimes_{\mathbb{F}}V)_{{\la}k\ra}=(U(P)(1\otimes_{\mathbb{F}}
V))_{{\la}k\ra}$
 if and only if  the vectors
 $$\{
 p_{i_{1}}\cdots
p_{i_{k}}(1\otimes v_{j}) \ | \ 1 \leq i_{1} \leq \cdots \leq i_{k}
\leq n, j \in \overline{1,\ell} \}\eqno(3.33)$$ are linearly
independent. Equivalently, the corresponding  $\ell\ell_{k}\times
\ell \ell_{k}$ system of homogeneous linear equations
$$M_{k}X=0\eqno(3.34)$$ has only zero solution by (3.31).

Suppose $$\sum\limits_{1 \leq i_{1}\leq i_{2}\leq \cdots \leq
i_{k}\leq n, j \in \overline{1,d}}
a_{i_{1},i_{2},\cdots,i_{k}}^{j}p_{i_{1}}p_{i_{2}} \cdots p_{i_{k}}
(1\otimes v_{j})=0, \ a_{i_{1}, i_{2}, \cdots , i_{k}}^{j} \in
\mathbb{F}.\eqno(3.35)$$ It can be written as
$$\sum\limits_{1 \leq i_{1}\leq i_{2}\leq \cdots \leq i_{k}\leq n}
p_{i_{1}}p_{i_{2}}\cdots p_{i_{k}}(1\otimes w_{i_{1}, i_{2}, \cdots
, i_{k}})=0, \ w_{i_{1}, i_{2}, \cdots , i_{k}} \in V.\eqno(3.36)$$
 Then $$\partial_{x_{l}}.\sum\limits_{1 \leq i_{1}\leq
i_{2}\leq \cdots \leq i_{k}\leq n} p_{i_{1}}p_{i_{2}}\cdots
p_{i_{k}}(1\otimes w_{i_{1}, i_{2}, \cdots , i_{k}})=0, \  \forall \
l \in \overline{1,n}.\eqno(3.37)$$ Equivalently,
$$
\sum\limits_{1 \leq i_{1}\leq i_{2}\leq \cdots \leq i_{k}\leq
n}\sum\limits_{s=1}^{n}p_{i_{1}}p_{i_{2}}\cdots\hat{p}_{i_{s}}\cdots
p_{i_{k}}\triangle_{l,i_{s}}^{k}(1\otimes w_{i_{1}, i_{2}, \cdots ,
i_{k}})=0\eqno(3.38)$$ by (3.20). Moreover, it can be written as the
form
$$ \sum\limits_{1 \leq j_{1}\leq j_{2}\leq \cdots \leq j_{k-1}\leq
n} p_{j_{1}}p_{j_{2}}\cdots p_{j_{k-1}}(1\otimes u_{j_{1},
j_{2}, \cdots , j_{k-1}})=0,\eqno(3.39)$$where
\begin{eqnarray*}\hspace{2cm}& &
1\otimes u_{j_{1}, j_{2}, \cdots ,
j_{k-1}}\\
&=&\sum\limits_{i=1}^{j_{1}}\triangle_{l,i}^{k}(1\otimes w_{i,j_{1},
j_{2}, \cdots ,
j_{k-1}})+\sum\limits_{i=j_{1}}^{j_{2}}\triangle_{l,i}^{k}(1\otimes
w_{j_{1},i,j_{2}, \cdots , j_{k-1}})+\cdots\\&&+
\sum\limits_{i=j_{s-1}}^{j_{s}}\triangle_{l,i}^{k}(1\otimes
w_{j_{1},j_{2}, \cdots ,j_{s-1},i,j_{s},\cdots j_{k-1}})+\cdots\\
& &+\sum\limits_{i=j_{k-1}}^{n}\triangle_{l,i}^{k}(1\otimes
w_{j_{1},j_{2}, \cdots ,
j_{k-1},i})=N_{k-1}\varpi_{n}\hspace{5cm}(3.40)
\end{eqnarray*}
 with
$$ N_{k-1}={\small \left[\begin{array}{cccc}
\sum\limits_{i=1}^{n}x_{i}\partial_{x_{i}}+x_{1}\partial_{x_{1}}+k-1&x_{2}\partial_{x_{1}}
&\cdots&x_{n}\partial_{x_{1}}\\
x_{1}\partial_{x_{2}}&\sum\limits_{i=1}^{n}x_{i}\partial_{x_{i}}+x_{2}
\partial_{x_{2}}+k-1&
\cdots&x_{n}\partial_{x_{2}}\\
\vdots&\vdots&\cdots&\vdots\\
x_{1}\partial_{x_{n}}&x_{2}\partial_{x_{n}}&\cdots&\sum\limits_{i=1}^{n}
x_{i}\partial_{x_{i}}+x_{n}
\partial_{x_{n}}+k-1
\end{array}
\right]}\eqno(3.41)$$ and
$$
\varpi_{n}= \left[\begin{array}{c} X_{1,j_{1}-1}\\
X_{j_{1},j_{2}-1}\\X_{j_{2},j_{3-1}}\\\vdots\\
X_{j_{k-1},n}\end{array} \right], \eqno(3.42)$$ in which
$$
X_{1,j_{1}-1}= \left[\begin{array}{c} 1 \otimes w_{1,j_{1},j_{2},
\cdots , j_{k-1}} \\ 1\otimes w_{2,j_{1},j_{2}, \cdots ,
j_{k-1}}\\\vdots \\1\otimes w_{j_{1}-1,j_{1}, j_{2},\cdots ,
j_{k-1}}\end{array} \right], \;\;X_{j_{1},j_{2}-1}=
\left[\begin{array}{c} 2\otimes w_{j_{1},j_{1},j_{2}, \cdots ,
j_{k-1}}\\1\otimes w_{j_{1},j_{1}+1,j_{2}, \cdots , j_{k-1}}
\\\vdots \\1\otimes w_{j_{1},j_{2}-1,j_{2}, \cdots , j_{k-1}}
\end{array} \right],\eqno(3.43) $$
$$X_{j_{2},j_{3-1}}=\left[\begin{array}{c} 2\otimes
w_{j_{1},j_{2},j_{2}, \cdots , j_{k-1}}
 \\1\otimes w_{j_{1},j_{2},j_{2}+1, \cdots , j_{k-1}}
\\\vdots \\1\otimes w_{j_{1},j_{2},j_{3}-1, \cdots , j_{k-1}}
\end{array} \right],\cdots,\;\;
 X_{j_{k-1},n}=\left[\begin{array}{c}
2\otimes w_{j_{1},j_{2}, \cdots , j_{k-1},j_{k-1}}\\1\otimes
w_{j_{1},j_{2}, \cdots ,j_{k-1}, j_{k-1}+1}\\\vdots \\1\otimes
w_{j_{1},j_{2}, \cdots , j_{k-1},n}
\end{array}
\right]. \eqno(3.44)$$ We know that $N_{k-1}\varpi_{n}=0$ has only
zero solution if and only if $\mu_{i}+|\mu|-i+k \neq 0,  \ \forall \
i \in I_{1}$. By Lemma 3.2.1 and the claim, the Lemma is followed.
\hspace{1cm} $\Box$\vspace{0.3cm}

Next, we study the structure $(U(P)(1\otimes_{\mathbb{F}}
V))_{{\la}j\ra}$ as an $\bar L_n$-module and finally get a
necessary and sufficient condition for the irreducibility of  $L_{n+1}$-module  ${\cal A}\otimes_{\mathbb{F}}V$ (cf. Proposition 3.2.5 ).\vspace{0.3cm}

 Note that $\bar L_n \simeq gl(n)$  ( resp. $\bar L_n'\simeq sl(n)$ ) according to (3.5). Obviously, (3.12) implies that
 as $\overline{L}_{n}$-module
$({\cal A}\otimes_{\mathbb{F}}V)_{{\la}j\ra} $ is isomorphic  to
tensor product module $V(j\epsilon_{1})\otimes_{\mathbb{F}}V$
$$({\cal
A}\otimes_{\mathbb{F}}V)_{{\la}j\ra} \cong
V(j\varepsilon_{1})\otimes_{\mathbb{F}}V(\mu)=\bigoplus_{\underline{c}\in
I(\mu,j)} V(\mu+\underline{c})\eqno(3.45)$$ (cf. Lemma 2.1.3).
Denote
$$ p^{\underline{c}}=p_{1}^{c_1}p_{2}^{c_2}\cdots p_{n}^{c_n}\qquad \mbox{for} \
\underline{c}=(c_{1},\cdots,c_{n}) \in \mathbb{N}^{n}.\eqno(3.46)$$
Define the linear map
$$\varphi_{j}:({\cal A}\otimes_{\mathbb{F}}V)_{{\la}j\ra}\rightarrow (U(P)(1\otimes_{\mathbb{F}} V))_{{\la}j\ra}$$$$\sum\limits_{|\underline{l}|=j, \  q \in
\overline{1,l}}a_{\underline{l},q}x^{\underline{l}}\otimes v_{q} \mapsto \sum\limits_{|\underline{l}|=j, \  q \in
\overline{1,l}}a_{\underline{l},q}p^{\underline{l}}.(1\otimes  v_{q}). \eqno(3.47)$$It is easy to verify that
$\varphi_{j}$ is an $\overline{L}_{n} $-module homomorphism by (3.21).\vspace{0.3cm}

For any $\underline{c}\in I(\mu,j)$ (cf. (2.11)), let
$\xi_{\underline{c}}$ be a maximal vector
 for  highest weight module $V(\mu+\underline{c})\subseteq ({\cal A}\otimes_{\mathbb{F}}V)_{{\la}j\ra}$.
 Then
$$\hat{\xi}_{\underline{c}}=\varphi_{j}(\xi_{\underline{c}})\eqno(3.48)$$
is also a maximal vector of $V(\mu+\underline{c}) \subseteq
(U(P)(1\otimes_{\mathbb{F}} V))_{{\la}j\ra}$.
 Write
$$\hat{\xi}_{\underline{c}}=q_{\underline{c}}\xi_{\underline{c}},\; \ \ q_{\underline{c}}\in \mathbb{F}.\eqno(3.49)$$We know that
$$(U(P)(1\otimes_{\mathbb{F}} V))_{{\la}j\ra}\bigcap
V(\mu+\underline{c}) \neq \{0\}\Rightarrow
V(\mu+\underline{c})\subseteq (U(P)(1\otimes_{\mathbb{F}}
V))_{{\la}j\ra}.\eqno(3.50)$$
So
$$(U(P)(1\otimes_{\mathbb{F}} V))_{{\la}j\ra}\bigcap
V(\mu+\underline{c}) \neq \{0\} \qquad \mbox{iff} \qquad  q_{\underline{c}} \neq 0.\eqno(3.51)$$In the following Lemma, we will calculate $q_{\underline{c}} $ explicitly:
 \vspace{0.3cm}

{\bf Lemma 3.2.3 } {\it  For any $\underline{c}\in I(\mu,j)$, we have $$q_{\underline{c}}=\prod\limits_{s=1}^{n}\prod\limits_{i=1}^{c_{s}}(\mu_{s}+|\mu|-s+i),
\eqno(3.52)$$where we treat
$\prod\limits_{i=1}^{c_{s}}(\mu_{s}+|\mu|-s+i)=1$ if $c_{s}=0$. }
\vspace{0.3cm}

{\it Proof}\quad  Let $v_{\mu}$ be a highest weight vector of the
given $gl(n)$-module $V$.

For any $\underline{c}\in I(\mu,j)$,  we write a maximal vector
$\xi_{\underline{c}}$
 for  highest weight module $V(\mu+\underline{c})$ as
$$\xi_{\underline{c}}=a_{\underline{c}}x^{{\underline{c}}}\otimes v_{\mu}+
\sum\limits_{(\nu,\underline{l})\neq (\mu, \underline{c} ),
{\underline{c}}+\mu={\underline{l}}+\nu, \ i \in
\overline{1,m({\nu})}}a_{\nu,\underline{l}}^{i}x^{\underline{l}}\otimes
v_{\nu}^{i},\qquad 0 \neq a_{\underline{c}},  \ a_{\nu,\underline{l}}^{i}
\in \mathbb{F},\eqno(3.53)$$where  $m({\nu})$ denotes the
multiplicity of weight $\nu \in \Pi_\mu$.\vspace{0.3cm}

{\it Case 1.  $
c_{n}\neq 0$.}\vspace{0.3cm}

Note that $[x_{s}\partial_{x_{t}},\partial_{x_{n}}]=0, \ \forall \ 1 \leq s
< t \leq n$. We know
$$0 \neq
\partial_{x_{n}}. \xi_{\underline{c}}=
a_{\underline{c}}c_{n}x^{{\underline{c}}-\varepsilon_{n}}\otimes v_{\mu}+
\sum\limits_{(\nu,\underline{l})\neq (\mu, \underline{c} ),
{\underline{c}}+\mu={\underline{l}}+\nu, \ i \in
\overline{1,m({\nu})}}a_{\nu,\underline{l}}^{i}l_{n}x^{\underline{l}
-\varepsilon_{n}}\otimes v_{\nu}^{i}\eqno(3.54)$$ is a maximal
vector for the  $\overline{L}_{n}$- highest weight module
$V(\mu+\underline{c}-\varepsilon_{n})$. Set $$\xi_{\underline{c}-\varepsilon_{n}}=\partial_{x_{n}}. \xi_{\underline{c}}.\eqno(3.55)$$Obviously,
$$\hat{\xi}_{\underline{c}}=a_{\underline{c}}p^{{\underline{c}}}.(1\otimes v_{\mu})+
\sum\limits_{(\nu,\underline{l})\neq (\mu, \underline{c} ),
{\underline{c}}+\mu={\underline{l}}+\nu, \ i \in
\overline{1,m({\nu})}}a_{\nu,\underline{l}}^{i}p^{\underline{l}}.(1\otimes
v_{\nu}^{i}),\eqno(3.56)$$
$$
\hat{\xi}_{\underline{c}-\varepsilon_{n}}=\varphi_{j}(\partial_{x_{n}}. \xi_{\underline{c}})
=a_{\underline{c}}c_{n}p^{
{\underline{c}}-\varepsilon_{n}}.(1\otimes v_{\mu})+
\sum\limits_{(\nu,\underline{l})\neq (\mu, \underline{c} ),
{\underline{c}}+\mu={\underline{l}}+\nu, \ i \in
\overline{1,m({\nu})}}a_{\nu,\underline{l}}^{i}l_{n}p^{
\underline{l}-\varepsilon_{n}}.(1\otimes v_{\nu}^{i}).\eqno(3.57)
$$
Write $$ \partial_{x_{n}}.
\hat{\xi}_{\underline{c}}=b_{\underline{c}}(n)
\hat{\xi}_{\underline{c}-\varepsilon_{n}}, \qquad b_{\underline{c}}(n) \in \mathbb{F}.\eqno(3.58)$$By
(3.21) and (3.54), we have
\begin{eqnarray*}& &
\partial_{x_{n}}.
\hat{\xi}_{\underline{c}}=\partial_{x_{n}}.\varphi_{j}( \xi_{\underline{c}})\\
&=&c_{n}
a_{\underline{c}}p^{{\underline{c}}-\varepsilon_{n}}\Delta_{n,n}^{j}(
1 \otimes v_{\mu})\\&&+ \sum\limits_{(\nu,\underline{l})\neq (\mu,
\underline{c} ), {\underline{c}}+\mu={\underline{l}}+\nu, \ i \in
\overline{1,m({\nu})}}a_{\nu,\underline{l}}^{i}[\sum\limits_{s=1}
^{n-1}l_{s}p^{\underline{l}-\varepsilon_{s}}(1\otimes
E_{s,n}v_{\nu}^{i})+l_{n}p^{\underline{l}-\varepsilon_{n}}\Delta_{n,n}^{j}(1\otimes
v_{\nu}^{i})] .\hspace{1cm}(3.59)
\end{eqnarray*}
 Denote
$$E_{s,n}.v_{\mu+\varepsilon_{n}-\varepsilon_{s}}^{i}=\Im_{i}v_{\mu}, \qquad \Im_{i} \in \mathbb{F}.\eqno(3.60)$$
Then, we have
$$b_{\underline{c}}(n)=\frac{c_{n}a_{\underline{c}}(j-1+\mu_{n}+|\mu|)
+\sum\limits_{s=1}^{n-1}(1+c_{s})\sum\limits_{i=1}^{m(\mu-\varepsilon_{s}
+\varepsilon_{n})}
a_{\mu-\varepsilon_{s}+\varepsilon_{n},\underline{c}+\varepsilon_{s}-\varepsilon_{n}}^{i}
\Im_{i}}{c_{n}a_{\underline{c}}}\eqno(3.61)$$ by (3.58)-(3.61). For any $s
\in \overline{1,n-1}$, we get
\begin{eqnarray*}& &
0=x_{s}\partial_{x_{n}}.
\xi_{\underline{c}}\\
&=&c_{n}a_{\underline{c}}x^{\underline{c}+\varepsilon_{s}-\varepsilon_{n}}\otimes
v_{\mu}+ \sum\limits_{(\nu,\underline{l})\neq (\mu, \underline{c} ),
{\underline{c}}+\mu={\underline{l}}+\nu, \ i \in
\overline{1,m({\nu})}}a_{\nu,\underline{l}}^{i}[l_{n}x^
{\underline{l}+\varepsilon_{s}-\varepsilon_{n}}\otimes
v_{\nu}^{i}+x^{\underline{l}}\otimes
E_{s,n}.v_{\nu}^{i}].\hspace{0.2cm}(3.62)
\end{eqnarray*}
Therefore,
$$c_{n}a_{\underline{c}}+\sum\limits_{i=1}^{m(\mu-\varepsilon_{s}+\varepsilon_{n})}
a_{\mu-\varepsilon_{s}+\varepsilon_{n},\underline{c}+\varepsilon_{s}-
\varepsilon_{n}}^{i}\Im_{i}=0\eqno(3.63)$$ Hence,  we have
$$b_{\underline{c}}(n)=c_{n}-n+\mu_{n}+|\mu|.\eqno(3.64)$$Furthermore,
by (3.55), we obtain
$$  \partial_{x_{n}}.
\hat{\xi}_{\underline{c}}=q_{\underline{c}}\partial_{x_{n}}.\xi_{\underline{c}}
=q_{\underline{c}}\xi_{\underline{c}-\varepsilon_{n}}=
b_{\underline{k}}(n)
\hat{\xi}_{\underline{c}-\varepsilon_{n}}
=b_{\underline{k}}(n)q_{\underline{c}-\varepsilon_{n}}\xi_{\underline{c}-\varepsilon_{n}}
\eqno(3.65)$$which implies that
$$q_{\underline{c}}=b_{\underline{c}}(n)q_{\underline{c}-\varepsilon_{n}}\eqno(3.66)$$

{\it Case 2. $c_{n}=0$.}

Suppose $c_{n-1}\neq 0$. We claim that  $l_{n}=0$ for any
$\underline{l}$ of $ a_{\nu,\underline{l}}^{i} \neq 0$ in (3.54).
Indeed, we can write $$ (\nu_{1}-\nu_{2},\cdots,\nu_{n-1}-\nu_{n})
=\sum\limits_{i=1}^{n-1}(\mu_{i}-\mu_{i+1})\varpi_{i}
-\sum\limits_{i=1}^{n-1}k_i\alpha_i.\eqno(3.67)$$ Since $
\nu_{n}+l_{n}=\mu_{n}+c_{n}$, we have
$l_{n}=\mu_{n}-\nu_{n}=\mu_{n}-(\mu_{n}+k_{n-1})$ by (2.9). Thus we
get $l_{n}+k_{n-1}=0$. So $l_{n}=0$.

Hence,
$$x_{n-1}\partial_{x_{n}}.\partial_{x_{n-1}}.\xi_{\underline{c}}
=-\partial_{x_{n}}.\xi_{\underline{c}}+\partial_{x_{n-1}}.
x_{n-1}\partial_{x_{n}}.\xi_{\underline{c}}=0.\eqno(3.68)$$This
implies that  $\partial_{x_{n-1}}.\xi_{\underline{c}}$ is  a maximal
vector for the highest weight $\overline{L}_{n}$-module
$V(\mu+\underline{c}-\varepsilon_{n-1})$.  Repeating the process
from (3.56) to (3.64), we get
$$q_{\underline{c}}=(\mu_{n-1}+|\mu|-(n-1)+c_{n-1})
q_{\underline{c}-\varepsilon_{n-1}}.\eqno(3.69)$$ Then induction
implies that $$ q_{\underline{c}}=
\prod\limits_{s=1}^{n}\prod\limits_{i=1}^{c_{s}}(\mu_{s}+|\mu|-s+i).
\eqno(3.70)$$Thus the lemma is proved. \hspace{1cm}
$\Box$\vspace{0.3cm}

 By Lemma 3.2.3, we know that $(U(P)(1\otimes_{\mathbb{F}} V))_{{\la}j\ra}= ({\cal
A}\otimes_{\mathbb{F}}V)_{{\la}j\ra}$ iff  $q_{\underline{c}}\neq 0$ for any $\underline{c}\in I(\mu,j), \ 0 < j \in \mathbb{N}$. Thus we  get the following sufficient and necessary condition for the irreducibility of $L_{n+1}$-module  ${\cal A}\otimes_{\mathbb{F}}V$:
 \vspace{0.3cm}

{\bf Proposition 3.2.4 } {\it The $L_{n+1}$-module  ${\cal A}\otimes_{\mathbb{F}}V$ is irreducible iff  $\ q_{\underline{c}}\neq 0$ for any $\underline{c}\in I(\mu,j), \ 0 < j \in \mathbb{N}$.
}\vspace{0.3cm}

\subsection{Jordan-Holder Series for
  $L_{n+1}$-module ${\cal A}\otimes_{\mathbb{F}}V$}

\qquad In this section, we will assume that ${\cal A}\otimes_{\mathbb{F}}V \neq U(P)(1\otimes_{\mathbb{F}}
V)
 $ and prove the irreducibility of the quotient module $({\cal A}\otimes_{\mathbb{F}}V) /
U(P)(1\otimes_{\mathbb{F}} V)$.\vspace{0.3cm}

First, we study the structure of $({\cal A}\otimes_{\mathbb{F}}V) /
U(P)(1\otimes_{\mathbb{F}} V)$ as an $\overline{L}_{n} $-module (cf.
Lemma 3.3.1 and Lemma 3.3.2). Recall the $\overline{L}_{n} $-module
homomorphism $\varphi_{j}$ defined by (3.47). Denote
$\mbox{Ker}(\varphi_{j})$ by ${\cal{R}}_{{\la}j\ra}$ and set
 $$I(\mu,j)'=\{\underline{c}\in I(\mu,j) \ | \ \ q_{\underline{c}} =0\}.\eqno(3.71)$$According to Lemma 3.2.3, we have the following result:\vspace{0.3cm}

{\bf Lemma 3.3.1} {\it \quad For  any $  1 \leq j \in
\mathbb{N}$, we have ${\cal{R}}_{{\la}j\ra}=\bigoplus\limits_{\underline{c}\in I(\mu,j)'}V(\mu+\underline{c})$ and the following direct sum of submodules for
$\overline{L}_{n} $:
$ ({\cal A}\otimes_{\mathbb{F}}V)_{{\la}j\ra}= (U(P)(1\otimes_{\mathbb{F}} V))_{{\la}j\ra}\oplus
{\cal{R}}_{{\la}j\ra}.$
}
\vspace{0.3cm}

Let $1 \leq k \in \mathbb{N}$. Suppose that
 $$({\cal A}\otimes_{\mathbb{F}}V)_{{\la}i\ra}=(U(P)(1\otimes V))_{{\la}i\ra} \;\;  \mbox{for any
} \ \;i \in \overline{0,k},\eqno(3.72)$$ but $$({\cal
A}\otimes_{\mathbb{F}}V)_{{\la}j\ra}\neq (U(P)(1\otimes
V))_{{\la}j\ra} \;\;\mbox{when}\ j \geq k+1. \eqno(3.73)$$ Then by
Lemma 3.3.1, as $\overline{L}_{n} $-modules,
$$({\cal A}\otimes_{\mathbb{F}}V) /
U(P)(1\otimes_{\mathbb{F}} V) \cong\bigoplus\limits_{j=
k+1}^\infty{\cal{R}}_{{\la}j\ra}.\eqno(3.74)$$

{\bf  Lemma 3.3.2} {\it \quad Let $1 \leq k \in \mathbb{N}$ such
that (3.72) and (3.73) hold. Then $\bar{L}_n$-module
${\cal{R}}_{{\la}k+1\ra}$ is irreducible.}\vspace{0.3cm}

{\it Proof}\quad Since $({\cal
A}\otimes_{\mathbb{F}}V)_{{\la}k\ra}=(U(P)(1\otimes V))_{{\la}k\ra}
$, Lemma 3.2.3 implies
$$ q_{\underline{c}}\neq 0 \qquad  \forall \  \underline{c}\in I(\mu,k)
.\eqno(3.75)$$
For convenience, we denote
$$q_{s}(\underline{c})=\prod\limits_{i=1}^
{c_{s}}(\mu_{s}+|\mu|-s+i).\eqno(3.76)$$So
$$q_{\underline{c}}=\prod\limits_{s=1}^ {n}q_{s}(\underline{c}).\eqno(3.77)$$
Assume $\underline{m}\in I(\mu,k+1)'$, i.e.
$q_{\underline{m}}=0$. Since the $gl(n)$-modules
$$
V((k+1)\varepsilon_{1})\otimes_{\mathbb{F}}V\simeq ({\cal A}\otimes_{\mathbb{F}}V)_{{\la}k+1\ra}\subset
V(\varepsilon_1)\otimes_{\mathbb{F}} ({\cal
A}\otimes_{\mathbb{F}}V)_{{\la}k\ra}\simeq V(\varepsilon_{1})\otimes_{\mathbb{F}}(V
(k\varepsilon_{1})
\otimes_{\mathbb{F}}V) \eqno(3.78) $$ in the sense of
monomorphism,
 we know that $$\exists \quad \underline{t}\in {\mathbb{N}}^{n} \
\mbox{and} \ r \in \overline{1,n}, \  \mbox{such \ that}  \
\underline{t}\in I(\mu,k) \ \mbox{and} \
\underline{m}=\underline{t}+\varepsilon_{r}.\eqno(3.79)$$ Therefore,
$m_{r}=t_{r}+1$. The fact $q_{\underline{m}}=0$ implies that
$q_{r}(\underline{m})=0$ because
$q_{s}(\underline{m})=q_{s}(\underline{t})\neq 0$ for $s \neq r$.
Then $q_{r}(\underline{t})\neq 0$ and
$q_{r}(\underline{m})=0=q_{r}(\underline{t})(\mu_{r}+|\mu|-r+t_{r}+1)$
imply that
$$\mu_{r}+|\mu|-r+t_{r}+1=0=\mu_{r}+|\mu|-r+m_{r}.\eqno(3.80)$$
Obviously, $1 \leq m_{r}=t_{r}+1 \leq |\underline{m}|=k+1$ and $\underline{m}\in I(\mu,k+1)'\subset I(\mu,k+1) $ imply that $$ m_{r}\leq
\mu_{r-1}-\mu_{r}\eqno(3.81)$$ by (2.11).
Assume $1 \leq m_{r}< k$. Then by (2.11) and (3.81), we know
there exists some $\underline{l}\in {\mathbb{N}}^{n}$ satisfying
$l_{r}=m_{r}$ and $\underline{l}\in I(\mu,k) $. So
$q_{r}(\underline{l})=q_{r}(\underline{m})=0$. Furthermore,
$q_{\underline{l}}=0$, which contradicts (3.75). Therefore,
$m_{r}=k+1$, i.e. $\underline{m}=(k+1)\epsilon_{r}$.

Suppose there exists another $ (k+1)\varepsilon_{s}\in I(\mu,k+1)'$ but $s \neq r$. Then
$0=\mu_{r}+|\mu|-r+k+1=\mu_{s}+|\mu|-s+k+1$. This is impossible, since $\mu_{s}\geq \mu_{r}$ whenever $s< r$. Thus we prove that $|I(\mu,k+1)'|=1$, i.e.
$\bar{L}_n$-module ${\cal{R}}_{{\la}k+1\ra}$ is irreducible.\hspace{1cm} $\Box$\vspace{0.3cm}

Set$$I_{s}=\{1\}\bigcup\{j \in \overline{2,n} \ | \
\mu_{j-1}-\mu_{j} \geq s\}.\eqno(3.82)$$
By the above two lemmas, we can  give the proof of (i) in the Main Theorem:\vspace{0.3cm}

{\bf Proposition 3.3.3} \quad {\it The vector space ${\cal
A}\otimes_{\mathbb{F}}V$ is an irreducible $sl(n+1)$-module if and
only if
$$\forall \  1 \leq s \in \mathbb{N}, \ \mu_{i}+|\mu|-i+s \neq 0
\qquad\mbox{for\ any}\;\;  i\in I_{s}.\eqno(3.83)$$} \vspace{0.1cm}

{\it Proof }
Assume that there exist $1 \leq s \in \mathbb{N} $ and $i\in I_{s}$
satisfying
 $\mu_{i}+|\mu|-i+s  = 0 $. Then $V(\mu+s\varepsilon_{i})\subseteq {\cal{R}}_{{\la}s\ra}$ by (2.11), (3.71), (3.82) and
 Lemma 3.3.1. Hence, $({\cal
A}\otimes_{\mathbb{F}}V)_{{\la}s\ra}\neq (U(P)(1\otimes
V))_{{\la}s\ra} $, i.e. $sl(n+1)$-module ${\cal
A}\otimes_{\mathbb{F}}V$ is reducible.

Suppose that $sl(n+1)$-module ${\cal A}\otimes_{\mathbb{F}}V$ is
reducible. Let $k$ satisfying (3.72) and (3.73). By Lemma 3.3.2, we
have $V(\mu+(k+1)\epsilon_{s})={\cal{R}}_{{\la}k+1\ra}$ for some $s
\in I_{k+1}$. So Lemma 3.2.3 implies $\mu_{s}+|\mu|-s+k+1 = 0
$.\hspace{1cm} $\Box$\vspace{0.3cm}

{\bf Remark 3.3.4}\quad The condition (3.83) is equivalent  to (1.7) and
(1.8) given in the Main Theorem.\vspace{0.3cm}

In the rest of this section, we study the relationship between any
$\bar{L}_n$-module ${\cal{R}}_{{\la}j\ra}$ and  $\bar{L}_n$-module
${\cal{R}}_{{\la}j+1\ra}$  for any $j\geq k+1$ based on the
decomposition of tensor  module and the projection operator
techniques for $gl(n)$ (Recall Lemma 2.1.3 and projection operators
appeared in Section 2.2). And we finally prove the irreducibility of
quotient module $({\cal A}\otimes_{\mathbb{F}}V) /
U(P)(1\otimes_{\mathbb{F}} V)$ (cf. Proposition 3.3.8)
.\vspace{0.3cm}

{\bf  Lemma 3.3.5} {\it \quad Let $k+1\leq s\in\mathbb{N}$. For
$\nu=\mu+\underline{l}$ with $\underline{l} \in I(\mu,s+1)'$, there
exist $ \nu'=\mu+\underline{m}$ with  $\underline{m} \in I(\mu,s)'$
and $r \in \overline{1,n}$ such that $\nu=\nu'+\varepsilon_{r}$ and
$\mu_{r}+|\mu|-r+m_{r}+1\neq 0$. }\vspace{0.3cm}

{\it Proof}\quad Suppose $\nu=\mu+\underline{l}$ with
$\underline{l}\in I(\mu,s+1)'$. Then
$q_{\underline{l}}=0$.\vspace{0.1cm}

 {\it Claim. There exists $r \in  \overline{1,n}$ such that $\mu_{r}+|\mu|-r+l_{r}\neq 0$ and $l_{r}\geq 1$.}\vspace{0.1cm}

Set $$I_{\underline{l}}=\{t \in \overline{1,n} \ | \
\mu_{t}+|\mu|-t+l_{t}=0\}.\eqno(3.84)$$ It follows that
$|I_{\underline{l}}|=0,1$. Otherwise,
$\mu_{t}+|\mu|-t+l_{t}=0=\mu_{q}+|\mu|-q+l_{q}$ for some $t< q$.
This is impossible because $\mu+\underline{l}$ is a highest weight
implies that $l_{t}+\mu_{t}\geq l_{q}+\mu_{q}$ whenever $t< q$. It
is obvious that the claim holds when $|I_{\underline{l}}|=0$.

Now assume $|I_{\underline{l}}|=1$ and $r_{0}\in I_{\underline{l}}$. If the claim does not hold, then
$\underline{l}=(s+1)\varepsilon_{r_{0}}$ and $\mu_{t}+|\mu|-t \neq 0$ for any $t \neq r_{0}$.
 From Lemma 3.3.2, we know ${\cal{R}}_{{\la}k+1\ra}=V(\mu+(k+1)\varepsilon_{t})$ for some $t \in I_{k+1}$.
 Thus, Lemma 3.2.3 implies $q_{(k+1)\epsilon_{t}}=\prod\limits_{i=1}^{k+1}(\mu_{t}+|\mu|-t+i)=0$. Assume
  $\mu_{t}+|\mu|-t+r=0$ for some $r \in \overline{1,k+1}$. On the other hand,
  $\mu_{r_{0}}+|\mu|-r_{0}+s+1=0$ by (3.84) due to $r_{0}\in I_{\underline{l}}$. Hence, we have $\mu_{r_{0}}-\mu_{t}=r_{0}+r-(t+s+1)$. If $r_{0} \leq t$,
  then $\mu_{r_{0}}-\mu_{t} \geq 0$; which contradicts $r_{0}+r-(t+s+1)= r_{0}-t+r-(s+1)< 0$. If $r_{0} > t$, then
  $\mu_{t}+l_{t}=\mu_{t} \geq \mu_{r_{0}}+l_{r_{0}}=\mu_{r_{0}}+s+1$ because $\mu+\underline{l}=\mu+(s+1)\varepsilon_{r_{0}}$ is
  a highest weight. Therefore, $\mu_{t}-\mu_{r_{0}}\geq s+1$; i.e. $t+s+1-(r_{0}+r)\geq s+1$. Hence, $r_{0}-t+r \leq 0$. A contradiction arises. Thus  the claim holds.

Suppose that $r$ satisfies the claim. Take $\nu'=\nu-\varepsilon_{r}, \
\underline{m}=\underline{l}-\varepsilon_{r}$.
 We claim that $\underline{m}\in I(\mu,s)'$. In fact,  $\underline{m}\in I(\mu,s)$ by (2.11)
 because $l_{r}-1\leq s$ and $l_{r}\leq \mu_{r-1}-\mu_{r}$ implies $ l_{r}-1 < \mu_{r-1}-\mu_{r}$. Furthermore,
 $q_{\underline{l}}=(\mu_{r}+|\mu|-r+l_{r})q_{\underline{m}}=0$ implies that $q_{\underline{m}}=0$. Therefore, $\underline{m}\in I(\mu,s)'$. Thus the Lemma follows.
\hspace{1cm} $\Box$\vspace{0.3cm}

{\bf Lemma 3.3.6} {\it \quad We have
$\partial_{l}({\cal{R}}_{{\la}j\ra})\not\subseteq (U(P)(1\otimes
V))_{{\la}j-1\ra}$ for any $l \in \overline{1,n}$ and $j >
k+1.$}\vspace{0.3cm}

{\it Proof}\quad Assume
 $$\partial_{l}({\cal{R}}_{{\la}j\ra})\subseteq (U(P)(1\otimes_{\mathbb{F}} V))_{{\la}j-1\ra}\ \mbox{for \ some}   \ l \in
 \overline{1,n}\;\mbox{and}\;
  j > k+1.\eqno(3.85)$$
By (3.11), we know that
$$\partial_{l}(U(P)(1\otimes_{\mathbb{F}} V)_{{\la}j\ra})\subseteq
(U(P)(1\otimes_{\mathbb{F}} V))_{{\la}j-1\ra}.\eqno(3.86)$$ Then (3.85), (3.86) and Lemma 3.3.1 imply that
 $$\partial_{l}(({\cal A}\otimes_{\mathbb{F}}V)_{{\la}j\ra})\subseteq (U(P)(1\otimes
 V))_{{\la}j-1\ra},\eqno(3.87)$$ that is,
$$\partial_{l}(({\cal A}\otimes_{\mathbb{F}}V)_{{\la}j\ra})
\varsubsetneq ({\cal
A}\otimes_{\mathbb{F}}V)_{{\la}j-1\ra}.\eqno(3.88)$$
This contradicts the fact
$$\partial_{l}(({\cal
A}\otimes_{\mathbb{F}}V)_{{\la}j\ra}) = ({\cal
A}\otimes_{\mathbb{F}}V)_{{\la}j-1\ra}\;\; \ \mbox{for any} \ l \in
\overline{1,n}. \eqno(3.89)$$Thus the lemma follows.\hspace{1cm}
$\Box$\vspace{0.3cm}

Let $\{e_{1}',\cdots,e_{n}'\}$ be a basis for
$\overline{L}_{n}$-module $V(\varepsilon_{1})$. For any $1 \leq j
\in \mathbb{N}$, we define the following linear map:
$$T_{j} :V(\varepsilon_{1})\otimes_{\mathbb{F}}({\cal
A}\otimes_{\mathbb{F}}V)_{{\la}j\ra}\rightarrow({\cal
A}\otimes_{\mathbb{F}}V)_{{\la}j+1\ra}$$$$T_{j}(e_{i}'\otimes
v)=p_{i}.v, \ \forall \  v \in ({\cal
A}\otimes_{\mathbb{F}}V)_{{\la}j\ra}.\eqno(3.90)$$

{\bf Lemma 3.3.7 } \quad {\it The linear map  $\  T_{j} \ $  is an intertwining  operator from the tensor module $V(\varepsilon_{1})\otimes_{\mathbb{F}}({\cal
A}\otimes_{\mathbb{F}}V)_{{\la}j\ra}$ to  $({\cal
A}\otimes_{\mathbb{F}}V)_{{\la}j+1\ra}$ for $\overline{L}_{n}$. }\vspace{0.3cm}

{\it Proof}\quad For any $s,t,i \in \overline{1,n}$ and $v \in
({\cal A}\otimes_{\mathbb{F}}V)_{{\la}j\ra}$, we have
\begin{eqnarray*}& & T_{j} (x_{s}\partial_{x_{t}}.(e_{i}'\otimes v))\\&=&T_{j} (x_{s}\partial_{x_{t}}(e_{i}')\otimes v+e_{i}'\otimes x_{s}\partial_{x_{t}}.v)=T_{j}(\delta_{t,i}e_{s}'\otimes v+e_{i}'\otimes x_{s}\partial_{x_{t}}.v)\\&=&\delta_{t,i}p_{s}.v+p_{i}.x_{s}\partial_{x_{t}}.v
=x_{s}\partial_{x_{t}}.p_{i}.v=x_{s}\partial_{x_{t}}.T_{j}(e_{i}'\otimes
v)\hspace{5.1cm}(3.91)
\end{eqnarray*}
by (3.21). Thus the lemma follows.\hspace{1cm} $\Box$\vspace{0.3cm}

Based on the Lemma 3.3.5, Lemma 3.3.6 and Lemma 3.3.7, we can prove (ii) of the Main Theorem in the following:\vspace{0.3cm}

{\bf Proposition 3.3.8 } {\it If ${\cal A}\otimes_{\mathbb{F}}V \neq
U(P)(1\otimes_{\mathbb{F} } V)
 $, then
 $\{0\} \subset U(P)(1\otimes_{\mathbb{F}} V) \subset {\cal A}\otimes_{\mathbb{F}}V$ is a Jordan-Holder Series for
  $L_{n+1}$-module ${\cal A}\otimes_{\mathbb{F}}V$.}\vspace{0.3cm}

 {\it Proof} \quad
Suppose  that $W+U(P)(1\otimes_{\mathbb{F}} V)$ is any nonzero
submodule of quotient module $({\cal A}\otimes_{\mathbb{F}}V) /
U(P)(1\otimes_{\mathbb{F}} V)$, where $W=\bigoplus\limits_{j\geq
k+1}W\bigcap {\cal{R}}_{{\la}j\ra}$ is a weighted subspace of ${\cal
A}\otimes_{\mathbb{F}}V$. By Lemma 3.3.6, we get $W\bigcap
{\cal{R}}_{{\la}k+1\ra}\neq \{0\}$.

By Lemma 3.3.5, we know that for any $s \geq k+1$ and
$\nu=\mu+\underline{l}$ with $\underline{l} \in I(\mu,s+1)'$, there
exist $ \nu'=\mu+\underline{m}$ with $\underline{m} \in I(\mu,s)'$
and $r \in \overline{1,n}$ such that $\nu=\nu'+\varepsilon_{r}$ and
$\mu_{r}+|\mu|-r+m_{r}+1\neq 0$. Therefore, the highest weight
module $V(\nu)\; (\subseteq {\cal{R}}_{{\la}s+1\ra} )$ of highest
weight $\nu$ appears in the decomposition of
$\overline{L}_{n}$-tensor module
$V(\varepsilon_{1})\otimes_{\mathbb{F}}V(\nu')$ ($ \subseteq
V(\varepsilon_{1})\otimes_{\mathbb{F}}{\cal{R}}_{{\la}s\ra}
\subseteq V(\varepsilon_{1})\otimes_{\mathbb{F}}({\cal
A}\otimes_{\mathbb{F}}V)_{{\la}s\ra} $). \vspace{0.3cm}

{\it Claim. \quad  There exists some maximal vector $v_{\nu}$ (resp.
$\xi_{\underline{m}+\varepsilon_{r}}$ ) of highest weight module
$V(\nu)\subseteq V(\varepsilon_{1})\otimes_{\mathbb{F}}V(\nu') $
(resp. $V(\nu) \subseteq {\cal{R}}_{{\la}s+1\ra} $ ) satisfying
$$T_{s}(v_{\nu})=(\mu_{r}+|\mu|-r+m_{r}+1)\xi_{\underline{m}+\varepsilon_{r}}\neq
0.\eqno(3.92)$$} Assume that $w_{\nu}$ is a maximal vector of
irreducible module $V(\nu)\subseteq
V(\varepsilon_{1})\otimes_{\mathbb{F}}V(\nu')$. Since
$T_{s}:V(\varepsilon_{1})\otimes_{\mathbb{F}}({\cal
A}\otimes_{\mathbb{F}}V)_{{\la}s\ra}\rightarrow({\cal
A}\otimes_{\mathbb{F}}V)_{{\la}s+1\ra}$ is an intertwining  operator
for $\overline{L}_{n}$, we know  $T_{s}(w_{\nu})$ is also a maximal
vector of $\overline{L}_{n}$-module $V(\nu)\subseteq
{\cal{R}}_{{\la}s+1\ra} \;( \subseteq ({\cal
A}\otimes_{\mathbb{F}}V)_{{\la}s+1\ra})$. Since the  maximal vector
$w_{\nu}$ must take the following form:
$$w_{\nu}=a_{r}e_{r}'\otimes \varpi_{\nu'}+
\sum_{j=1}^{r-1}\sum_{i=1}^{m(\nu'-\varepsilon_{j}+\varepsilon_{r})}
a_{j}^{i}e_{j}'\otimes
v_{\nu'-\varepsilon_{j}+\varepsilon_{r}}^{i},\qquad 0 \neq a_{r}, \
a_{j}^{i} \in \mathbb{F},\eqno(3.93)$$ where $\varpi_{\nu'}$ is a
maximal vector of the highest weight module $V(\nu')$ and the set
$\{v_{\nu'-\varepsilon_{j}+\varepsilon_{r}}^{i}\mid
i\in\ol{1,m(\nu'-\varepsilon_{j}+\varepsilon_{r}))}\}$ is a basis
 of the weight subspace
 $[V(\nu')]_{\nu'-\varepsilon_{j}+\varepsilon_{r}}$.
 For
convenience, we write
$$v_{\nu'-\varepsilon_{j}+\varepsilon_{r}}^{i}=\sum\limits_{ \varepsilon_{j}-\varepsilon_{r}=\sum\limits_{k=1}^{m}i_{k}
(\varepsilon_{s_{k}}-\varepsilon_{t_{k}});\: s_{i}> t_{i}
}a_{\underline{s},\underline{t}}^{\underline{i}}(x_{s_{1}}
\partial_{t_{1}})^{i_{1}}.\cdots.(x_{s_{m}}\partial_{t_{m}})^{i_{m}}.
\varpi_{\nu'}.\eqno(3.94)
$$
Therefore,$$T_{s}(w_{\nu})=a_{r}p_{r}.\varpi_{\nu'}+
\sum_{j=1}^{r-1}\sum_{i=1}^{m(\nu'-\varepsilon_{j}+\varepsilon_{r})}a_{j}^{i}p_{j}.
v_{\nu'-\varepsilon_{j}+\varepsilon_{r}}^{i}.\eqno(3.95)$$So (3.94) and (3.95) imply that
$$T_{s}(w_{\nu})=n_{\varepsilon_{r}}.\varpi_{\nu'} \qquad  \mbox{for \ some} \  n_{\varepsilon_{r}} \in U(P\oplus \bar{L}_{n}).\eqno(3.96)$$
Let $\varpi_{\nu'}={\xi}_{\underline{m}}$ (resp.
$\varpi_{\nu'}=\hat{\xi}_{\underline{m}}$ ) in (3.96), where
$$\xi_{\underline{m}}=a_{\underline{m}}x^{{\underline{m}}}\otimes
v_{\mu}+ \sum\limits_{(\eta,\underline{l})\neq (\mu, \underline{m}
), {\underline{m}}+\mu={\underline{l}}+\eta, \ i \in
\overline{1,m({\eta})}}a_{\eta,\underline{l}}^{i}x^{\underline{l}}\otimes
v_{\eta}^{i},\qquad 0 \neq a_{\underline{m}},  \
a_{\eta,\underline{l}}^{i} \in \mathbb{F};\eqno(3.97)$$
$$\hat{\xi}_{\underline{m}}=a_{\underline{m}}p^{{\underline{m}}}.(1\otimes v_{\mu})+
\sum\limits_{(\eta,\underline{l})\neq (\mu, \underline{m} ),
{\underline{m}}+\mu={\underline{l}}+\eta, \ i \in
\overline{1,m({\eta})}}a_{\eta,\underline{l}}^{i}p^{\underline{l}}.(1\otimes
v_{\eta}^{i}),\qquad 0 \neq a_{\underline{m}},  \ a_{\eta,\underline{l}}^{i}
\in \mathbb{F}.\eqno(3.98)$$
Take $$\hat{\xi}_{\underline{m}+\varepsilon_{r}}=n_{\varepsilon_{r}}.
\hat{\xi}_{\underline{m}}.\eqno(3.99)$$
Then by Lemma 3.2.3, we have$$\hat{\xi}_{\underline{m}+\varepsilon_{r}}=q_{\underline{m}+\varepsilon_{r}}
{\xi}_{\underline{m}+\varepsilon_{r}}=
n_{\varepsilon_{r}}.
\hat{\xi}_{\underline{m}}=q_{\underline{m}}n_{\varepsilon_{r}}.
{\xi}_{\underline{m}}.\eqno(3.100)$$Hence, we have
$$(\mu_{r}+|\mu|-r+m_{r}+1){\xi}_{\underline{m}+\varepsilon_{r}}=
n_{\varepsilon_{r}}. {\xi}_{\underline{m}}.\eqno(3.101)$$ Now we
take
\begin{eqnarray*}v_{\nu}&=&
\sum_{j=1}^{r-1}\sum_{i=1}^{m(\nu'-\varepsilon_{j}+\varepsilon_{r})}\sum\limits_{
\varepsilon_{j}-\varepsilon_{r}=\sum\limits_{k=1}^{m}i_{k}
(\varepsilon_{s_{k}}-\varepsilon_{t_{k}});\: s_{i}> t_{i}
}a_{j}^{i}a_{\underline{s},\underline{t}}^{\underline{i}}e_{j}'\otimes(x_{s_{1}}
\partial_{t_{1}})^{i_{1}}.\cdots.(x_{s_{m}}\partial_{t_{m}})^{i_{m}}.
{\xi}_{\underline{m}}\\ &
&+a_{r}e_{r}'\otimes
{\xi}_{\underline{m}}.\hspace{11.2cm}(3.102)\end{eqnarray*}
Then $ T_{s}(v_{\nu})=(\mu_{r}
+|\mu|-r+m_{r}+1)\xi_{\underline{m}+\varepsilon_{r}}\neq 0$ by
(3.101) and (3.102). Therefore, (2.30) and (3.92) imply that
$$\{0\}\neq
(T_{s}|_{V(\varepsilon_{1})\bigotimes_{\mathbb{F}}V(\nu')}\circ
\tilde{P}_{r})(V(\varepsilon_{1})\otimes_{\mathbb{F}}V(\nu'))=V(\nu),\eqno(3.103)$$
where $$\tilde{P}_{r}=\prod\limits_{l\neq
r}(\frac{\tilde{M}-\tilde{d}_{l}}{\tilde{d}_{r}-\tilde{d}_{l}}), \
\tilde{d}_{i}=i-1-m_{i}-\mu_{i},\eqno(3.104)$$ and  $\tilde{M}$ is
the matrix in (2.27). Hence, ${\cal{R}}_{{\la}s\ra} \subseteq W$ for
any $s \geq k+1$.
 Thus we prove  $({\cal A}\otimes_{\mathbb{F}}V) / U(P)(1\otimes_{\mathbb{F}} V)$ is
irreducible.\hspace{1cm} $\Box$ \psp

Suppose ${\cal{R}}_{{\la}k+1\ra}=V(\mu+(k+1)\varepsilon_{r})$ for
some $r \in I_{k+1}$. From the proof of Proposition 3.3.8, we know $({\cal
A}\otimes_{\mathbb{F}}V) / U(P)(1\otimes_{\mathbb{F}} V)$ is
generated by
$\bar{\xi}_{\mu+(k+1)\varepsilon_{r}}={\xi}_{\mu+(k+1)
\varepsilon_{r}}+U(P)(1\otimes_{\mathbb{F}}
V)$, i.e. $({\cal A}\otimes_{\mathbb{F}}V) /
U(P)(1\otimes_{\mathbb{F}}
V)=U(L_{n+1}).\bar{\xi}_{\mu+(k+1)\varepsilon_{r}}$. Hence, we get
the following result:\vspace{0.3cm}

{\bf Proposition 3.3.9}\quad {\it The  $L_{n+1}$-module $({\cal
A}\otimes_{\mathbb{F}}V) / U(P)(1\otimes_{\mathbb{F}} V)$  is
isomorphic to the irreducible module $ U(P)(1\otimes_{\mathbb{F}}
M)$, where $M$ is a $gl(n)$-irreducible module admitting the
character $\chi_{\mu+(k+1)\varepsilon_{r}}$.}\vspace{0.3cm}

{\it Proof} \quad Since $\partial_{i}.
{\xi}_{\mu+(k+1)\varepsilon_{r}} \in (U(P)(1\otimes
V))_{{\la}k\ra}=({\cal A}\otimes_{\mathbb{F}}V)_{{\la}k\ra} $, we
get $\partial_{i}. \bar{\xi}_{\mu+(k+1)\varepsilon_{r}}=0$ for any
$i \in \overline{1,n}$. It follows from Proposition 3.3.8 that both
$({\cal A}\otimes_{\mathbb{F}}V) / U(P)(1\otimes_{\mathbb{F}} V)$
and $ U(P)(1\otimes_{\mathbb{F}} M)$ are irreducible
$L_{n+1}$-modules. So it is easy to verify that
$$\vartheta:({\cal A}\otimes_{\mathbb{F}}V) / U(P)(1\otimes_{\mathbb{F}} V)\rightarrow U(P)(1\otimes_{\mathbb{F}} M);
x.\bar{\xi}_{\mu+(k+1)\varepsilon_{r}}\mapsto x.(1 \otimes
v_{\mu+(k+1)\varepsilon_{r}})\eqno(3.105)$$ is an $L_{n+1}$-module
isomorphism, where $x \in U(P\oplus\bar{L}_{n})$ and
$v_{\mu+(k+1)\varepsilon_{r}}$ is a maximal vector of
$M$.\hspace{1cm} $\Box$ \psp

From Proposition 3.3.3, Proposition 3.3.8 and Proposition 3.3.9, we get the Main Theorem.\vspace{0.3cm}

{\bf Remark 3.3.10}\quad The irreducible module $U(P)(1\otimes_{\mathbb{F}} V)$ is cyclic, i.e. it is generated by one vector. From Lemma 3.1.1 and (3.12), we know that $U(P)(1\otimes_{\mathbb{F}} V)$ is in general not a highest weight module. We can easily verify the following result from the Main Theorem:\vspace{0.3cm}

{\bf Corollary 3.3.11}\quad {\it Assume one of the conditions in (1.7) and (1.8) fails. Denote $$i_{0}=\mbox{min}\{i\in \overline{1,n} \ | \ \mu_{i}+|\mu|-i+1 \in -\mathbb{N}\}.\eqno(3.106)$$ We have:

(i) The integer $k=-\mu_{i_{0}}-|\mu|+ i_{0}-1$ satisfies  (3.72) and (3.73).

(ii) The irreducible module $U(P)(1\otimes_{\mathbb{F}} V)$ is finite dimensional highest weight module with highest weight $k\omega_{1}+\sum\limits_{i=2}^{n}m_{i-1}\omega_{i}$ iff  $ \ i_{0}=1$, where $m_{i}=\mu_{i}-\mu_{i+1}$ for $i \in \overline{1,n-1}$.

}
\end{document}